\documentclass[final,1p,times,numbers,sort&compress]{elsarticle}
\usepackage{txfonts}

\usepackage{amssymb}
\usepackage{amsthm}
\usepackage{amscd}
\usepackage{amsmath}
\usepackage{amsfonts}
\usepackage{amssymb}
\usepackage{graphicx}

\numberwithin{equation}{section}

\newcommand{\ra}{\rightarrow}

\newcommand{\f}{\frac}

\newcommand{\be}{\begin{equation}}
\renewcommand{\ra}{\rightarrow}
\newcommand{\ee}{\end{equation}}
\newcommand{\bea}{\begin{eqnarray}}
\newcommand{\eea}{\end{eqnarray}}
\newcommand{\bna}{\begin{eqnarray*}}
\newcommand{\ena}{\end{eqnarray*}}

\renewcommand{\le}{\left}
\newcommand{\ri}{\right}

\newcommand{\ep}{\epsilon}

\journal{$\ast\ast\ast$}

\begin{document}

\begin{frontmatter}

\title{Existence of solutions for a higher order Kirchhoff type problem with exponential critical growth}

\author[label1,label2]{Liang Zhao\corref{cor1}}
 \ead[label1]{ liangzhao@bnu.edu.cn}
\address[label1]{School of Mathematical Sciences,
 Beijing Normal
University, Beijing 100875, P. R. China}
\address[label2]{Laboratory of Mathematics and
Complex Systems, Ministry of Education,
Beijing 100875, P. R. China}
\cortext[cor1]{Corresponding author.}
\author[label3]{Ning Zhang}
 \ead[label1]{nzhang@amss.ac.cn}
\address[label3]{China Institute for Actuarial Science,
 Central University of Finance and Economics, Beijing 100081, P. R. China}

\begin{abstract}

The higher order Kirchhoff type equation
$$\int_{\mathbb{R}^{2m}}(|\nabla^m u|^2
+\sum_{\gamma=0}^{m-1}a_{\gamma}(x)|\nabla^{\gamma}u|^2)dx
\le((-\Delta)^m u+\sum_{\gamma=0}^{m-1}(-1)^\gamma
\nabla^\gamma\cdot(a_\gamma (x)\nabla^\gamma u)\ri)
=\frac{f(x,u)}{|x|^\beta}+\ep h(x)\ \  \text{in}\ \  \mathbb{R}^{2m}$$
is considered in this paper. We assume that the nonlinearity of the equation has exponential critical growth and prove that, for a positive $\ep$ which is small enough, there are two distinct nontrivial solutions to the equation. When $\ep=0$, we also prove that the equation has a nontrivial mountain-pass type solution.
\end{abstract}


\end{frontmatter}

\section{Introduction and main results}
Let
$\nabla^\gamma u$, $\gamma\in \{0,1,2,\cdots,m\}$, be the $\gamma$-th order gradients of a function $u\in W^{m,2}(\mathbb{R}^{2m})$ which are defined by
$$\nabla^\gamma u:=\le\{\begin{array}{ll}
   \Delta^{\f{\gamma}{2}}u & \gamma\ \  \text{even},\\[1.5ex]
   \nabla\Delta^{\f{\gamma-1}{2}}u & \gamma\ \  \text{odd}.
   \end{array}\ri.$$
Here and throughout this paper, $m\geq 2$ is an even integer and we use the notations that
$$\Delta^0 u=\nabla^0 u=u.$$

Consider the following nonlinear functional
\be\label{functional}
J_\ep(u) =\f{1}{4}\le(\int_{\mathbb{R}^{2m}}(|\nabla^m u|^2
+\sum_{\gamma=0}^{m-1}a_{\gamma}(x)|\nabla^{\gamma}u|^2)dx\ri)^2
-\int_{\mathbb{R}^{2m}}\frac{F(x,u)}{|x|^\beta}dx
-\ep\int_{\mathbb{R}^{2m}}h u dx
\ee
which is related to a higher order nonlocal partial differential equation
\be\label{equation}
\int_{\mathbb{R}^{2m}}(|\nabla^m u|^2
+\sum_{\gamma=0}^{m-1}a_{\gamma}(x)|\nabla^{\gamma}u|^2)dx
\le((-\Delta)^m u+\sum_{\gamma=0}^{m-1}(-1)^\gamma
\nabla^\gamma\cdot(a_\gamma (x)\nabla^\gamma u)\ri)
=\frac{f(x,u)}{|x|^\beta}+\ep h.
\ee
Here $\ep$ is a nonnegative constant, $h(x)\nequiv 0$ belongs to the dual space of $E$ which will be defined later, $0\leq \beta<2m$ and $a_\gamma(x)$ are continuous functions satisfying\vspace{.2cm}

\noindent${\bf (A_1)}$there exist positive constants $a_\gamma$, $\gamma=0,1,2,\cdots, m-1$, such that
$a_\gamma(x)\geq a_\gamma$ for all $x\in \mathbb{R}^{2m}$;

\noindent${\bf (A_2)}$ $(a_0(x))^{-1}\in L^1(\mathbb{R}^{2m})$.\vspace{.2cm}

Since the equation contains an integral over $\mathbb{R}^{2m}$, it is no longer a pointwise identity and should be dealt with as a nonlocal problem. We call (\ref{equation}) a higher order Kirchhoff type equation because it is related to the stationary analog of the equation
\begin{equation}\label{equation2}
\rho\frac{\partial^2 u}{\partial t^2}-\left(\frac{\rho_0}{h}+\frac{E}{2L}\int_0^L \left|\frac{\partial u}{\partial x}\right|^2 dx\right)\frac{\partial^2 u}{\partial x^2}=0,
\end{equation}
where $\rho$, $\rho_0$, $h$, $E$ and $L$ are constants. This equation was presented by Kirchhoff \cite{Kir} as an extension of the classical D'Alembert wave equation for free vibrations of elastic string produced by transverse vibrations. This kind of nonlocal problem also appears in other fields, for example, biological systems where $u$ describes a process which depends on the average of itself (for instance, population density). One can refer to \cite{AlvCor, AlvCorMa, PerZha} and the references therein for more details. After the work of Lions \cite{Lions}, where a functional analysis approach was proposed to this kind of equations, various models of Kirchhoff type have been studied by many authors using the variational framework, see, for example, \cite{Ane,CheKuoWu,Che,CheLi,CheWu,FanLiu,FigIkoJun,LiYan,LiWu,LiLiShi,
LiaZha,LiaLiShi,Nai,YeTan,ZhaSunNie,HeZou}, and the references therein. In particularly, Li and Yang \cite{LiYan} studied the following equation and proved the existence of at least two positive solutions.
$$\le\{\begin{array}{ll}
   M\left(\int_{\mathbb{R}^N}(|\nabla u|^N+V(x)|u|^N)dx\right)(-\Delta_N u+V(x)|u|^{N-2}u)=\lambda A(x)|u|^{p-2}u+f(u)\ \  x\in \mathbb{R}^N,\\[1.5ex]
   u\in W^{1,N}(\mathbb{R}^N),
   \end{array}\ri.$$
where $\Delta_N u= div (|\nabla u|^{N-2}\nabla u)$ is the $N$-Laplacian operator of $u$, $M(s)=s^k$ for $k>0$, $s\geq 0$, $1<p<N$, $\lambda>0$ is a real parameter, $A(x)$ is a positive function in $L^\sigma (\mathbb{R}^N)$ with $\sigma=\frac{N}{N-p}$, $V$ is a potential function and $f$ is a nonlinearity term having critical exponential growth.

On the other hand, similar variational methods are also used to study equations without the nonlocal integral. For example, see \cite{AdiYad,AdiYan,Cao,dooMedSev,FigdooRuf,FigMiyRuf,
GazGruSqu,LamLu1,Pan,ReiWet,Yan1,Yan2,YanZha,Zha1,Zha2,ZhaCha} and the reference therein. Among these results, the first author and Chang \cite{ZhaCha} proved an Adams type inequality and applied it to get the multiplicity result of a higher order quasilinear equation as following
\be\label{equation3}
(-\Delta)^m u+\sum_{\gamma=0}^{m-1}(-1)^\gamma\nabla^{\gamma}
\cdot(a_\gamma(x)\nabla^\gamma u)
=\f{f(x,u)}{|x|^\beta}+\epsilon h(x).
\ee
It is natural to ask the a question that can we generalize the results in \cite{LiYan} to higher order cases? This is the motivation of our paper. Our equation (\ref{equation}) is a higher order Kirchhoff type equation with a singular nonlinearity which is a nonlocal version of the equation (\ref{equation3}).
A primary tool to study this kind of equations is the Adams type inequality. Precisely, we need the following theorem which is proved in \cite{ZhaCha}.\vspace{.2cm}

\noindent{\bf Theorem A.} {\it Let $m\geq 2$ be an even integer and $0\leq\beta<2m$, then for any $0\leq\alpha\leq\le(1-\f{\beta}{2m}\ri)\alpha(m,2m)$,
\be\label{adams}
\sup_{u\in W^{m,2}(\mathbb{R}^{2m}),\|u\|_{E}\leq 1}
\int_{\mathbb{R}^{2m}}
\f{e^{\alpha u^2}-1}{|x|^\beta}dx<\infty,
\ee
where $\alpha(m,2m)=(4\pi)^m m!$. Furthermore, the inequality is sharp, which means that when $\alpha>\le(1-\f{\beta}{2m}\ri)\alpha(m,2m)$, the integrals are still finite for any $u\in E$, but the supremum goes to infinity.}\\

According to the variational structure of the functional (\ref{functional}), we assume that the nonlinearity $f(x,s):\mathbb{R}^{2m}\times\mathbb{R}\ra\mathbb{R}$ is a continuous function and satisfies the following growth conditions\vspace{.2cm}

\noindent${\bf (H_1)}$ There exist constants $\alpha_0$, $b_1$, $b_2>0$ and $\theta\geq 3$ such that for all $(x,s)\in \mathbb{R}^{2m}\times\mathbb{R}$,
$$|f(x,s)|\leq b_1|s|^3+b_2|s|^\theta (e^{\alpha_0 s^2}-1).$$
\vspace{.2cm}

\noindent${\bf (H_2)}$ There exists $\mu>4$ such that for all $x\in\mathbb{R}^{2m}$ and $s\neq 0$,
$$0<\mu F(x,s)\equiv\mu\int_0^s f(x,t)dt\leq sf(x,s).$$
\vspace{.2cm}


Define a function space
$$E:=\left\{u\in W^{m,2}(\mathbb{R}^{2m})| \int_{\mathbb{R}^{2m}}(|\nabla^m u|^2+\sum_{\gamma=0}^{m-1}a_{\gamma}(x)|\nabla^\gamma u|^2) dx<+\infty\right\}$$
and denote the norm of $u\in E$ by
$$\|u\|_E :=\left(\int_{\mathbb{R}^{2m}}(|\nabla^m u|^2+\sum_{\gamma=0}^{m-1}a_{\gamma}(x)|\nabla^\gamma u|^2) dx\right)^{1/2}.$$
Here and in the sequel we use $E^*$ to denote the dual space of $E$. Define
$$\lambda_p:=\inf_{u\in E\setminus \{0\}}\frac{\|u\|_E}{\left(\int_{\mathbb{R}^{2m}}\frac{u^p}{|x|^\beta} dx\right)^{1/p}}.$$
Obviously, we can conclude that $\lambda_p>0$ from the following proposition which can be found in \cite{ZhaCha}.\vspace{.2cm}

\noindent{\bf Proposition B.} {\it Under assumptions $(A_1)$ and $(A_2)$, we have that the space $E$ is compactly embedded into the space $L^q(\mathbb{R}^{2m})$ for any $q\geq 1$.}\\

We also assume\vspace{.2cm}

\noindent${\bf (H_3)}$
$\limsup_{s\ra 0}\f{4|F(x,s)|}{s^4}<\lambda_4^{4}$ uniformly with respect to $x\in \mathbb{R}^{2m}$.
\vspace{.2cm}

Let $\alpha(m,2m)=(4\pi)^m m!$ be the constant in Theorem A. We assume\vspace{.2cm}

\noindent${\bf (H_4)}$ There exist constants $p>4$ and $C_p$ such that
$$|f(s)|\geq C_p |s|^{p-1},$$
where
$$C_p>\le(\f{\mu(p-4)}{p(\mu-4)}\ri)^{\f{p-4}{4}}
\le(\f{\alpha_0}{\left(1-\frac{\beta}{2m}\right)\alpha(m,2m)}\ri)^{\f{p-4}{2}}\lambda_p^p.$$
\vspace{.2cm}

\noindent{\bf Remark.} To construct an example of $f(x,s)$ satisfying $(H_1)-(H_5)$, one can refer to examples in \cite{Zha1}.\vspace{.2cm}

We will see that the functional $J_\ep$ satisfies the geometric conditions of the mountain-pass theorem. Namely, there exist two constants $r_\ep>0$ and $\vartheta_\ep>0$ such that $J_\ep(u)\geq \vartheta_\ep$ when $\|u\|_E =r_\ep$. And there exists some $e\in E$ satisfying $\|e\|_E>r_\ep$ such that $J_\ep(e)<0$. Moreover, $J_\ep(0)=0$. The min-max level $C_M$ of $J_\ep$ is defined by
$$C_M=\min_{l\in \mathcal{L}}\max_{u\in l}J_\ep(u),$$
where $\mathcal{L}=\{l\in \mathcal{C}([0,1],E):l(0)=0,l(1)=e\}$. Thus we have the first result.\vspace{.2cm}

\noindent{\bf Theorem 1.1.} {\it Assume $(A_1)$, $(A_2)$ and $(H_1)-(H_4)$. There exists some $\ep_0>0$ such that, for $0\leq\ep<\ep_0$, the equation (\ref{equation}) has a nontrivial mountain-pass type weak solution $u_0$ and the min-max level $C_M$ has an upper bound
$$C_M< \le(\f{\mu-4}{4\mu}\ri)
\le(\f{\left(1-\frac{\beta}{2m}\right)\alpha(m,2m)}
{\alpha_0}\ri)^{2}.$$}\\

When $\ep=0$, (\ref{equation}) becomes
\be\label{equation4}
\|u\|_E^2
\le((-\Delta)^m u+\sum_{\gamma=0}^{m-1}(-1)^\gamma
\nabla^\gamma\cdot(a_\gamma (x)\nabla^\gamma u)\ri)
=\frac{f(x,u)}{|x|^\beta}
\ee
and the corresponding functional is
\be\label{functional4}
J(u) =\f{1}{4}\le(\int_{\mathbb{R}^{2m}}(|\nabla^m u|^2
+\sum_{\gamma=0}^{m-1}a_{\gamma}(x)|\nabla^{\gamma}u|^2)dx\ri)^2
-\int_{\mathbb{R}^{2m}}\frac{F(x,u)}{|x|^\beta}dx.
\ee
Theorem 1.1 states that when $\ep=0$, namely for the equation (\ref{equation4}), the mountain-pass solution still exists. But to find another nontrivial minimum type solution to (\ref{equation}) which is distinct from $u_0$, we need $\ep\neq 0$. Precisely, we have\vspace{.2cm}

\noindent{\bf Theorem 1.2.} {\it Assume $(A_1)$, $(A_2)$ and $(H_1)-(H_4)$. Then there exists $\ep_1>0$ such that, for any $0<\ep< \ep_1$, the equation (\ref{equation}) has at least two distinct weak solutions.}\\

We organize this paper as follows: In Section 2, we estimate the min-max level of functional (\ref{functional}) and prove Theorem 1.1. In section 3, we prove Theorem 1.2.

\section{Mountain-pass type solution}

First we claim that $J_\ep\in \textit{C}^1(E,\mathbb{R})$. Namely for $u_k\rightarrow u$ in $E$, we have $J_\ep(u_k)\rightarrow J_\ep(u)$ and $J_\ep '(u_k)\rightarrow J_\ep '(u)$ as $k\rightarrow \infty$. In fact, we have
\bna
\frac{1}{4}(\|u_k\|_E^4-\|u\|_E^4)
&=&\frac{1}{4}(\|u_k\|_E^2+\|u\|_E^2)
(\|u_k\|_E+\|u\|_E)
(\|u_k-u_0+u_0\|_E-\|u\|_E)\\
&\leq& C\|u_k-u_0\|_E\rightarrow 0.
\ena
and
$$
\int_{\mathbb{R}^{2m}}h (u_k-u)dx
\leq C\|h\|_{E^*}\|u_k-u_0\|_E\rightarrow 0.
$$
If we can prove
\be\label{c1}
\lim_{k\rightarrow\infty}\int_{\mathbb{R}^{2m}}
\frac{F(x,u_k)}{|x|^\beta}dx=\int_{\mathbb{R}^{2m}}
\frac{F(x,u)}{|x|^\beta}dx,
\ee
we will get $\lim_{k\rightarrow\infty}J_\ep(u_k)=J_\ep(u)$.

By $(H_1)$ and $(H_2)$, we have
\be\label{c2}
F(x,u_k)\leq C_1|u_k|^4+C_2|u_k|^{\theta+1}(e^{\alpha_0 u_k^2}-1).
\ee
On the other hand, by Proposition A, up to a subsequence, we can assume that $\lim_{k\rightarrow \infty} \|u_k-u\|_{L^q}=0$ for any $q\geq 1$. Since $0\leq\beta<2m$, for some fixed constant $K>0$, using H\"{o}lder's inequality, we have
\bea\label{c3}
\int_{\mathbb{R}^{2m}}
\frac{|u_k-u|^4}{|x|^\beta}dx&=&\int_{|x|<K}
\frac{|u_k-u|^4}{|x|^\beta}dx+\int_{|x|\geq K}
\frac{|u_k-u|^4}{|x|^\beta}dx\nonumber\\
&\leq& \|x^{-1}\|_{L^{\beta p}(|x|<K)}^\beta\|u_k-u\|_{L^{4q}(|x|<K)}^4
+K^{-\beta}\|u_k-u\|_{L^{4}}\rightarrow 0,
\eea
where $1/p+1/q=1$. Define a function $\phi(s):=s^{\theta+1}(e^{\alpha_0 s^2}-1)$. By the mean value theorem, we have
\bea\label{c4}
\left|\phi(|u_k|)-\phi(|u|)\right|\leq \left|\phi'(\xi)\right|\left|u_k-u\right|\leq \left(\phi'(|u_k|)+\phi'(|u|)\right)\left|u_k-u\right|
\eea
For $\alpha>0$, $r>1$ and any $r'>r$, it is easy to check that there exists a positive constant $C$ which only depends on $\alpha$, such that for all $s\in \mathbb{R}$,
\begin{equation}\label{mp4}
(e^{\alpha s^2-1})^p\leq C(e^{\alpha p' s^2}-1).
\end{equation}
In fact, this is a lemma in \cite{Zha1}. Then we have for any $r> 1$
\bna
\int_{\mathbb{R}^{2m}}(e^{\alpha_0 u_k^2}-1)^r dx&\leq&
\int_{\mathbb{R}^{2m}}(e^{2r'\alpha_0 (u_k-u)^2+2r'\alpha_0 u^2}-1) dx\\
&\leq& \frac{1}{p}\int_{\mathbb{R}^{2m}}(e^{2r'\alpha_0 (u_k-u)^2}-1) dx+\frac{1}{q}\int_{\mathbb{R}^{2m}}(e^{2r'\alpha_0 u^2}-1) dx,
\ena
where we have used the fact that, for $a,b\geq 0$, $1/p+1/q=1$, $ab-1\leq \frac{a^p-1}{p}+\frac{b^q-1}{q}$. Since $\lim_{k\rightarrow \infty}\|u_k-u\|_E=0$, Theorem A gives that
\be\label{c5}
\sup_{k}\int_{\mathbb{R}^{2m}}(e^{\alpha_0 u_k^2}-1)^r dx<+\infty.
\ee
Using (\ref{c4}), (\ref{c5}) and H\"{o}lder's inequality, we get
\be\label{c6}
\lim_{k\rightarrow\infty}\int_{\mathbb{R}^{2m}}
\frac{|u_k|^{\theta+1}(e^{\alpha_0 u_k^2}-1)}{|x|^\beta}dx=\int_{\mathbb{R}^{2m}}
\frac{|u|^{\theta+1}(e^{\alpha_0 u^2}-1)}{|x|^\beta}.
\ee
In view of (\ref{c2}), (\ref{c3}) and (\ref{c6}), the Lebesgue's dominated convergence theorem gives (\ref{c1}). The proof of $\lim_{k\rightarrow \infty}J_\ep '(u_k)=J_\ep '(u)$ is similar and we omit the details here.

To use the mountain-pass theory to discuss the existence of solutions to (\ref{equation}), first we prove the following lemma which ensures that the functional (\ref{functional}) satisfies the mountain-pass conditions.\vspace{.2cm}

\noindent{\bf Lemma 2.1.} {\it Under assumptions $(A_1)$, $(A_2)$ and $(H_1)-(H_3)$, we have that the functional (\ref{functional}) satisfies the geometric conditions of mountain-pass theorem, namely\\
(i) $J_\ep(0)=0$.\\
(ii) there exist constants $\ep_2$, $r_\ep$ and $\vartheta_\ep>0$ such that, for $0\leq \ep<\ep_2$,  $J_\ep(u)\geq \vartheta_\ep$ when $\|u\|_E =r_\ep$.\\
(iii) there exists some $e\in E$ satisfying $\|e\|_E>r_\ep$ such that $J_\ep(e)<0$.}\\

\noindent{\bf Proof.} Obviously (i) is true.

By $(H_3)$, there exist constants $\delta, \sigma>0$ such that, for any $|s|\leq \delta$ and $x\in \mathbb{R}^{2m}$
\begin{equation}\label{mp1}
|F(x,s)|\leq \left(\frac{\lambda_4^4}{4}-\frac{\sigma}{4}\right)|s|^4.
\end{equation}
By $(H_1)$ and $(H_2)$, there holds, for $|s|>\delta$
\begin{equation}\label{mp2}
|F(x,s)|\leq C |s|^q(e^{\alpha_0 s^2}-1),
\end{equation}
where $q>\theta+1\geq 4$ and $C$ is a constant depending on $b_1$, $b_2$, $\delta$ and $\theta$. It follows from (\ref{mp1}) and (\ref{mp2}) that, for any $s\in \mathbb{R}$ and $x\in \mathbb{R}^{2m}$
\begin{equation}\label{mp3}
|F(x,s)|\leq \left(\frac{\lambda_4^4}{4}-\frac{\sigma}{4}\right)|s|^4+C |s|^q(e^{\alpha_0 s^2}-1).
\end{equation}
By H\"{o}lder's inequality,Proposition A and (\ref{mp4}), we have
\bea\label{mp5}
\int_{\mathbb{R}^{2m}}\f{e^{\alpha_0 u^2}-1}{|x|^\beta}|u|^q dx
&\leq&\le(\int_{\mathbb{R}^{2m}}\f{(e^{\alpha_0 u^2}-1)^{p_1}}{|x|^{\beta p_1}}dx\ri)^{1/p_1}
\le(\int_{\mathbb{R}^{2m}}|u|^{qp_2}dx\ri)^{1/p_2}\nonumber\\
&\leq&C\le(\int_{\mathbb{R}^{2m}}\f{e^{\alpha_0p_1' u^2}-1}{|x|^{\beta p_1}}dx\ri)^{1/p_1}
\|u\|_E^{q},
\eea
where $1/p_1+1/p_2=1$, $p_1'>p_1$ and $0\leq \beta p_1<2m$. Let $\varrho^2=\le(1-\f{\beta p_1}{2m}\ri)\f{\alpha(m,2m)}{\alpha_0p_1'}$. Then Theorem A and (\ref{mp5}) imply that, for $\|u\|_E\leq \varrho$,
\be\label{mp6}
\int_{\mathbb{R}^{2m}}\f{e^{\alpha_0 u^2}-1}{|x|^\beta}|u|^q dx
\leq
C\|u\|_E^{q}.
\ee
On the other hand, we have
\begin{equation}\label{mp7}
\int_{\mathbb{R}^{2m}}\frac{|u|^4}{|x|^\beta}dx\leq \frac{\|u\|_E^4}{\lambda_4^4}.
\end{equation}
In view of (\ref{mp3}), (\ref{mp6}) and (\ref{mp7}), we obtain, for $\|u\|_E\leq \varrho$
$$J_\ep(u)\geq \frac{\|u\|_E^4}{4\lambda_4^4}-C\|u\|_E^q-\ep \|h\|_{E^*}\|u\|_E.$$
When $\ep=0$, obviously if we can take $\|u\|_E=r_0$, where $0<r_0<\left(\frac{1}{8C\lambda_4^4}\right)^{1/(q-4)}$, we have $J(u)=\vartheta_0>\frac{r^4}{8\lambda_4^4}>0$. When $\ep\neq 0$, if $\ep$ is sufficiently small, we can take $r_\ep=(16\ep \lambda_4^4\|h\|_{E^*})^{1/3}<r_0$ and $J_\ep(u)=\vartheta_\ep\geq J(u)-\ep \|h\|_{E^*}\|u\|_E>\ep r_\ep\|h\|_{E^*}$.

To prove (iii), first we claim that, for any $u\in E\setminus \{0\}$ and $B_R$ which is a ball with radius $R$, there exist a constant $C>0$ such that
\begin{equation}\label{mp8}
\int_{\mathbb{R}^{2m}}F(x,u)dx\geq C\int_{B_R}|u|^\mu dx+O(1).
\end{equation}
In fact, $(H_2)$ implies that, for any constant $s_0>0$ and $s>s_0$, we have
$$F(x,s)\geq \frac{s^\mu}{s_0^\mu}F(x,s_0).$$
Since $F(x,s)>0$ for $s\neq 0$ and $f(x,s)$ is a continuous function, there exist constants $C_1, C_2>0$ such that for any $(x,s)\in B_R\times \mathbb{R}$,
$$F(x,s)\geq C_1 |s|^\mu-C_2.$$
Then we have
$$\int_{\mathbb{R}^{2m}}F(x,u)dx\geq \int_{B_R}F(x,u)dx\geq C_1\int_{B_R}|u|^\mu dx-C_2|B_R|$$
and the claim is proved. On the other hand, take any $u\in E\setminus \{0\}$. We can find a constant $R_0$ such that $\int_{B_{R_0}} u^\mu dx=\Lambda_0>0$. Together this fact with (\ref{mp8}), we get
\bna
J_\ep(t u)&=&\frac{t^4}{4}\left(\int_{\mathbb{R}^{2m}}(|\nabla^m u|^2+\sum_{\gamma=0}^{m-1}a_{\gamma}(x)|\nabla^\gamma u|^2) dx\right)^2-\int_{\mathbb{R}^{2m}}\frac{F(x,t u)}{|x|^\beta}dx-\ep t\int_{\mathbb{R}^{2m}} h u dx\\
&\leq&\frac{t^4}{4}\|u\|_E^4-C_1 t^\mu\int_{B_{R_0}}\frac{u^\mu}{|x|^\beta} dx+\ep t\|h\|_{E^*}\|u\|_E+O(1)\\
&\leq&\frac{t^4}{4}\|u\|_E^4-C_1 t^u\frac{\Lambda_u}{R_0^\beta}+\ep t\|h\|_{E^*}\|u\|_E+O(1).
\ena
Since $\mu>4$, we have
\begin{equation}\label{mp9}
\lim_{t\rightarrow +\infty}J_\ep(t u)=-\infty,
\end{equation}
which gives (iii).
$\hfill\Box$\\

Next we estimate the mountain-pass level of the functional (\ref{functional}), which confirms part of the results in Theorem 1.1. \vspace{.2cm}

\noindent{\bf Lemma 2.2.} {\it Assume $(A_1)$, $(A_2)$ and $(H_1)-(H_4)$. There exists $\ep_3>0$ such that for $0\leq \ep<\ep_3$, the mountain-pass level $C_M$ of functional (\ref{functional}) satisfies
$$C_M< \le(\f{\mu-4}{4\mu}\ri)
\le(\f{\left(1-\frac{\beta}{2m}\right)\alpha(m,2m)}
{\alpha_0}\ri)^{2}.$$}\\

\noindent{\bf Proof.} We can choose a bounded sequence of functions $\{u_k\}\subset E$ such that
$$\int_{\mathbb{R}^{2m}}\f{|u_k|^p}{|x|^\beta}dx=1\quad{\rm and}\quad \|u_k\|_E\ra \lambda_p.$$
Then by Proposition A, we can assume that there exists a function $u_p$ such that
   $$\begin{array}{lll}
   u_k\rightharpoonup u_p &{\rm in}\quad E,\\[1.5ex]
   u_k\ra u_p &{\rm in}\quad L^q(\mathbb{R}^{2m})\quad{\rm for\ \ all}\quad q\in [1,+\infty),\\[1.5ex]
   u_k(x)\ra u_p(x) &{\rm a.e.\ \  in\ \ } \mathbb{R}^{2m}.
   \end{array}$$
   These imply that
   $$\int_{\mathbb{R}^{2m}}\f{|u_k|^p}{|x|^\beta}dx\ra \int_{\mathbb{R}^{2m}}\f{|u_p|^p}{|x|^\beta}dx=1.$$
   On the other hand, we have
   $$\|u_p\|_E\leq \liminf_{k\ra\infty}\|u_k\|_E=\lambda_p .$$
   Thus we get $\|u_p\|_E=\lambda_p.$ Define a function $M_\ep(t): [0, +\infty)\ra \mathbb{R}$ by
    $$M_\ep(t):=\f{t^4}{4}\left(\int_{\mathbb{R}^{2m}}(|\nabla^m u_p|^2+\sum_{\gamma=0}^{m-1}a_\gamma (x)|\nabla^\gamma u_p|^2)dx\right)^2
-\int_{\mathbb{R}^{2m}}\f{F(x,tu_p)}{|x|^\beta}dx
-\ep t\int_{\mathbb{R}^{2m}} h udx.$$
   By $(H_4)$ and $\int_{\mathbb{R}^{2m}}\f{|u_p|^p}{|x|^\beta}dx=1$, we have
   \bna
   M_\ep(t)
   &\leq&\f{t^4}{4}\left(\int_{\mathbb{R}^{2m}}(|\nabla^m u_p|^2+\sum_{\gamma=0}^{m-1}a_\gamma (x)|\nabla^\gamma u_p|^2)dx\right)^2
   -C_p\f{t^p}{p}\int_{\mathbb{R}^{2m}}\f{|u_p|^p}{|x|^\beta}dx
   +\ep t\|h\|_{E^*}\|u_p\|_E\\
   &=&\f{\lambda_p^4}{4}t^4-\f{C_p}{p}t^p
   +\ep\lambda_p\|h\|_{E^*}t\\
   &\leq&\f{(p-4)}{4p}\f{S_p^{4p/(p-4)}}{C_p^{4/(p-4)}}
   +\ep\lambda_p\|h\|_{E^*}t_0,
   \ena
where $t_0$ is a constant which belongs to $[0,+\infty)$ and is independent of the choice of $\ep$. From the definitions of $C_p$ and $S_p$ in $(H_4)$, by choosing $\ep$ small enough, we get the desired results immediately.
$\hfill\Box$\\

By Lemma 2.1 and 2.2, if $\ep$ is sufficiently small, the functional (\ref{functional}) satisfies the conditions of the mountain-pass theorem except for the Palais-Smale condition. The mountain-pass theorem without the Palais-Smale condition \cite{Rab} implies that we can find a Palais-Smale sequence $\{u_k\}\subset E$ at level $C_M$ and we have got an upper bound estimate for $C_M$. To get the existence result of a mountain-pass solution of (\ref{equation}), we only need to prove that there exists a function $u_0\in E$ such that $u_k\rightarrow u_0$ in $E$ as $k\rightarrow +\infty$. To this end, we need to estimate the norm of $u_k$ first. Precisely we have\vspace{.2cm}

\noindent{\bf Lemma 2.3.} {\it Assume $(A_1)$, $(A_2)$, and $(H_1)-(H_4)$.  Then for any Palais-Smale sequence $\{u_k\}\subset E$ of $J_\ep$ at level $C_M$, i.e.,
$$J_\ep(u_k)\ra C_M, J_\ep'(u_k)\ra 0\ \ \text{as}\ \ k\ra \infty,$$
there exists $\ep_4>0$ such that, for $0\leq \ep<\ep_4$, there holds
$$\limsup_{k\rightarrow \infty}\|u_k\|_E^2<\f{\left(1-\frac{\beta}{2m}\right)\alpha(m,2m)}
{\alpha_0}.$$}\\

\noindent{\bf Proof.} Since $\{u_k\}$ is a Palais-Smale sequence at level $C_M$, we have
\be\label{ps1}
\f{1}{4}\|u_k\|_E^4-\int_{\mathbb{R}^{2m}}\f{F(x,u_k)}{|x|^\beta}dx
-\ep\int_{\mathbb{R}^{2m}}h u_k dx\ra C_M\ \ \text{as}\ \ k\ra \infty,
\ee
and
\be\label{ps2}
\le|\|u_k\|_E^2\int_{\mathbb{R}^{2m}}(\nabla^{m}u_k\nabla^{m}\varphi
+\sum_{\gamma=0}^{m-1}a_\gamma(x)\nabla^\gamma u_k\nabla^\gamma\varphi)dx-\int_{\mathbb{R}^{2m}}\f{f(x,u_k)}{|x|^\beta}\varphi dx-\ep\int_{\mathbb{R}^{2m}}h \varphi dx
\ri|\leq \sigma_k\|\varphi\|_E,
\ee
where $\varphi$ is an arbitrary function in $C_0^\infty(\mathbb{R}^{2m})$ and $\sigma_k\ra 0$ as $k\ra\infty$. Multiplying (\ref{ps1}) by $\mu$ and and let $\varphi=u_k$ in(\ref{ps2}), we obtain
\begin{equation}\label{bd1}
\le(\f{\mu}{4}-1\ri)\|u_k\|_E^4
-\int_{\mathbb{R}^{2m}}\f{\mu F(x,u_k)-f(x,u_k)u_k}{|x|^\beta}dx\leq \mu C_M
+(\mu-1)\ep\|h\|_{E^*}\|u_k\|_E+o(\|u_k\|_E),
\end{equation}
where $\mu>4$ is the constant in $(H_2)$. Then (\ref{bd1}) and $(H_2)$ tell us that $u_k$ is bounded in $E$.

Furthermore, when $\ep=0$, we can get from (\ref{bd1}) and $(H_2)$ that
\begin{equation}\label{bd2}
\le(\f{\mu}{4}-1\ri)\|u_k\|_E^4\leq\mu C_M+o(1).
\end{equation}
Thus
$$
\limsup_{k\rightarrow +\infty}\|u_k\|_E^4\leq\frac{4\mu}{\mu-4} C_M.
$$
Using the estimate of $C_M$ in Lemma 2.2, we get
$$
\lim_{k\rightarrow +\infty}\|u_k\|_E^2<\f{\left(1-\frac{\beta}{2m}\right)\alpha(m,2m)}
{\alpha_0}.
$$

When $\ep\neq 0$, (\ref{bd2}) becomes
\begin{equation}\label{bd3}
\le(\f{\mu}{4}-1\ri)\|u_k\|_E^4\leq\mu C_M+(\mu-1)\ep\|h\|_{E^*}\|u_k\|_E+o(1).
\end{equation}
Applying Young's inequality, we have, for any $\eta>0$
$$
(\mu-1)\ep\|h\|_{E^*}\|u_k\|_E\leq \eta(\mu-1)\|u_k\|_E^4
+C_\eta\ep^{\frac{4}{3}}\|h\|_{E^*}^{\frac{4}{3}}
$$
where $C_\eta$ is a positive constant depending on $\eta$. We can conclude from Lemma 2.2 that there exists $\delta>0$ such that
\begin{equation}\label{bd4}
C_M=(1-\delta) \le(\f{\mu-4}{4\mu}\ri)
\le(\f{\left(1-\frac{\beta}{2m}\right)\alpha(m,2m)}
{\alpha_0}\ri)^{2}
\end{equation}
Choose $\eta=\frac{\delta(\mu-4)}{8(\mu-1)}$, (\ref{bd3}) gives
$$
\left(1-\frac{\delta}{2}\right)\le(\f{\mu}{4}-1\ri)\|u_k\|_E^4\leq\mu C_M+C_\eta\ep^{\frac{4}{3}}\|h\|_{E^*}^{\frac{4}{3}}+o(1).
$$
Therefore
$$
\limsup_{k\rightarrow \infty}\|u_k\|_E^4\leq\frac{8\mu C_M}{(2-\delta)(\mu-4)} +\frac{8C_\eta}{(2-\delta)(\mu-4)}\ep^{\frac{4}{3}}\|h\|_{E^*}^{\frac{4}{3}}.
$$
Using (\ref{bd4}) in the above inequality, we get
\begin{equation}\label{bd5}
\limsup_{k\rightarrow +\infty}\|u_k\|_E^4\leq
\le(\frac{2-2\delta}{2-\delta}\ri)
\le(\f{\left(1-\frac{\beta}{2m}\right)\alpha(m,2m)}
{\alpha_0}\ri)^{2}
+\frac{8C_\eta}{(2-\delta)(\mu-4)}\ep^{\frac{4}{3}}\|h\|_{E^*}^{\frac{4}{3}}.
\end{equation}
Now we choose $\ep_4$ such that it satisfies
$$
\ep_4^{\frac{4}{3}}=\frac{\delta(\mu-4)}{8C_\eta \|h\|_{E^*}^\frac{4}{3}}\le(\f{\left(1-\frac{\beta}{2m}\right)\alpha(m,2m)}
{\alpha_0}\ri)^{2}.
$$
It is easy to check that, for $0<\ep<\ep_4$, our estimate for $\|u_k\|_E$ still holds.
$\hfill\Box$\\

Next we prove that the functional $J_\ep$ satisfies the Palais-Smale condition at level $C_M$ and the proof of Theorem 1.1 can be completed by the following lemma.\vspace{.2cm}

\noindent{\bf Lemma 2.4.} {\it Assume $(A_1)$, $(A_2)$ and $(H_1)-(H_4)$. Let $\ep_0=\min\{\ep_i\}$, $i=2,3,4$. When $0\leq \ep<\ep_0$, we have, for any Palais-Smale sequence $\{u_k\}\subset E$ of $J_\ep$ at level $C_M\in \left(0,\le(\f{\mu-4}{4\mu}\ri)
\le(\f{\left(1-\frac{\beta}{2m}\right)\alpha(m,2m)}
{\alpha_0}\ri)^{2}\right)$, up to a subsequence, there exists $u_0\in E$ such that $u_k\rightarrow u_0$ in $E$ as $k\rightarrow +\infty$
and $u$ is a weak solution of (\ref{equation}).}\\

\noindent{\bf Proof.} By Lemma 2.3, we have,up to a subsequence,
\begin{equation}\label{ps1}
\lim_{k\rightarrow +\infty}\|u_k\|_E^2<\f{\left(1-\frac{\beta}{2m}\right)\alpha(m,2m)}
{\alpha_0}.
\end{equation}
Then Proposition A tells us that, up to a subsequence, there exists $u_0\in E$ such that
$$u_k\rightharpoonup u_0\ \ \text{in}\ \  E\ \ \text{and}\ \  u_k\ra u_0\ \ \text{in}\ \ L^q(\mathbb{R}^{2m})\ \ \text{for any}\ \ q\geq 1.$$
Since $\{u_k\}$ is a Palais-Smale sequence, we have, as $k\rightarrow 0$,
\bea\label{ps2}
\langle J_\ep '(u_k), u_k-u_0\rangle&=&\|u_k\|_E^2\int_{\mathbb{R}^{2m}}(\nabla^{m}u_k\nabla^{m}(u_k-u_0)
+\sum_{\gamma=0}^{m-1}a_\gamma(x)\nabla^\gamma u_k\nabla^\gamma(u_k-u_0))dx\nonumber\\
&&-\int_{\mathbb{R}^{2m}}\f{f(x,u_k)}{|x|^\beta}(u_k-u_0) dx-\ep\int_{\mathbb{R}^{2m}}h (u_k-u_0) dx\rightarrow 0.
\eea
Since $u_k\rightharpoonup u_0$ in $E$ and $\|u_k\|_E$ is bounded, we have, as $k\rightarrow 0$,
\be\label{ps3}
\|u_k\|_E^2\int_{\mathbb{R}^{2m}}(\nabla^{m}u_0\nabla^{m}(u_k-u_0)
+\sum_{\gamma=0}^{m-1}a_\gamma(x)\nabla^\gamma u_0\nabla^\gamma(u_k-u_0))dx\rightarrow 0.
\ee
Subtract (\ref{ps3}) from (\ref{ps2}), we get
\be\label{ps4}
\|u_k\|_E^2\|u_k-u_0\|_E^2-\int_{\mathbb{R}^{2m}}\f{f(x,u_k)}{|x|^\beta}(u_k-u_0) dx-\ep\int_{\mathbb{R}^{2m}}h (u_k-u_0) dx\rightarrow 0
\ee
Using H\"{o}lder's inequality, Theorem A, $(H_1)$, (\ref{mp4}) and (\ref{ps1}), similar to the proof of (\ref{mp5}), we have, for some $q>1$,
$$
\left|\int_{\mathbb{R}^{2m}}\f{f(x,u_k)}{|x|^\beta}(u_k-u_0) dx\right|\leq C\|u_k-u_0\|_{L^q}\rightarrow 0.
$$
On the other hand, $u_k\rightharpoonup u_0$ in $E$ gives
$$
\ep\int_{\mathbb{R}^{2m}}h (u_k-u_0) dx\rightarrow 0.
$$
Then we can deduce from (\ref{ps4}) that
$$
\|u_k\|_E^2\|u_k-u_0\|_E^2\rightarrow 0.
$$
Since the level $C_M>0$, we have $\lim_{k\rightarrow \infty}\|u_k\|\neq 0$. Therefore we know that $u_k\rightarrow u_0$ in $E$ and $u_0$ is a nontrivial solution of (\ref{equation}).
$\hfill\Box$\\

\section{Minimum type solution}

To get a solution different from $u_0$, first we estimate the infimum of $J_\ep$ near $0\in E$. We remark that in our proof, we can only deal with the case $\ep\neq 0$.\vspace{.2cm}

\noindent{\bf Lemma 3.1.} {\it Assume $(A_1)$, $(A_2)$ and $(H_2)$. Then there exists $\sigma>0$ such that
$$
\inf_{\|u\|_E\leq \sigma}J_\ep(u)=C_0<0.
$$}\\

\noindent{\bf Proof.} Since $E$ is a Hilbert space and $h\nequiv 0$, by the Riesz representation theorem, the equation
$$(-\Delta)^m u+\sum_{\gamma=0}^{m-1}(-1)^\gamma\nabla^{\gamma}
\cdot(a_\gamma(x)\nabla^\gamma u)
=\epsilon h(x).$$
has a nontrivial solution. Denote it by $v$ and we have
\be\label{jv1}
\|v\|_E^2=\int_{\mathbb{R}^{2m}}hvdx>0.
\ee
An easy computation shows that
$$\f{d}{dt}J_\ep(tv)=t^3\|v\|_E^4
-t^{-1}\int_{\mathbb{R}^{2m}}\f{f(x,tv)}{|x|^\beta}tvdx
-\ep\int_{\mathbb{R}^{2m}}hvdx.$$
$(H_2)$ gives that $f(x,tv)tv>0$ for any $t>0$. We have that there exists some $\eta>0$ such that for any $0<t<\eta$,
$$\f{d}{dt}J_\ep(tv)<0.$$
Noticing that $J_\ep(0)=0$, we have $J_\ep(tv)<0$ for any $0<t<\eta$. Take $\sigma=\eta\|v\|_E$, we get the lemma proved.
$\hfill\Box$\\

Define $B_{\sigma}=\{u\in E: \|u\|_E\leq \sigma\}.$ Suppose $\inf_{\|u\|_E\leq \sigma} J_\ep(u)=C_0$. By Lemma 3.1, we get $C_0<0$. In view of the facts that $\bar{B}_{\sigma}$ is a complete and convex metric space, $J_\ep$ is of class $\mathcal{C}^1$ and bounded from below, by Ekeland's variational principle, there exists a Palais-Smale sequence $\{\nu_k\}\subset \bar{B}_{\sigma}$ at level $C_0$. With out loss of generality, we can assume there exists $\nu_0$ such that
\be\label{converge}\begin{array}{lll}
   \nu_k\rightharpoonup \nu_0\ \ \text{in}\ \ E,\\[.5ex]
   \nu_k\ra \nu_0\ \ \text{in}\ \ L^q(\mathbb{R}^{2m})\ \ \text{for any}\ \ q\geq 1,\\[.5ex]
   \nu_k\ra \nu_0\ \ \text{a.e. in}\ \ \mathbb{R}^{2m}.
   \end{array}
\ee
Since $C_M>0$ and $C_0<0$, to complete the proof of Theorem 1.2, we only need to show the strongly convergence of $\nu_k$ to $\nu_0$ in $E$ as $k\rightarrow +\infty$.\vspace{.2cm}

\noindent{\bf Lemma 3.1.} {\it Assume $(A_1)$, $(A_2)$ and $(H_1)-(H_3)$. Then there exists $\ep_5>0$ such that, for $0<\ep<\ep_5$, the functional (\ref{functional}) satisfies the Palais-Smale condition, namely up to a subsequence, we have
$$
\|\nu_k-\nu_0\|_E\rightarrow 0\ \ \ \text{as}\ \ \ k\rightarrow +\infty.
$$}\\

\noindent{\bf Proof.} Similar to (\ref{bd3}), we have
\begin{equation}\label{mb1}
\le(\f{\mu}{4}-1\ri)\|\nu_k\|_E^4\leq\mu C_0+(\mu-1)\ep\|h\|_{E^*}\|\nu_k\|_E+o(1).
\end{equation}
Up to a subsequence, we have
\begin{equation}\label{mb2}
\le(\f{\mu}{4}-1\ri)\limsup_{k\rightarrow \infty}\|\nu_k\|_E^3-(\mu-1)\ep\|h\|_{E^*}< 0.
\end{equation}
Take
$$
\ep_5=\frac{\mu-4}{4\|h\|_{E^*}(\mu-1)}\left(\f{\left(1-\frac{\beta}{2m}\right)\alpha(m,2m)}
{\alpha_0}\right)^{\frac{3}{2}},
$$
for any $0<\ep<\ep_5$, we have
$$\limsup_{k\rightarrow \infty}\|\nu_k\|_E^2<\f{\left(1-\frac{\beta}{2m}\right)\alpha(m,2m)}
{\alpha_0}.$$
The same as the proof of Lemma 2.4, we can prove that $\lim_{k\rightarrow \infty}\|\nu_k-\nu_0\|_E=0$.
$\hfill\Box$\\

{\bf Acknowledgements.} This work is supported by NSFC 11201028 and NSFC 11426236. Part of this work was done during the first author visiting School of Mathematical Sciences, USTC. He thanks School of Mathematical Sciences, USTC for the hospitality and good working conditions.

\end{document}